\documentclass[11pt,a4paper]{article}
\usepackage[cp1251]{inputenc}
\usepackage{amsfonts, amssymb, amsthm, amsmath}
\usepackage{mathtext}
\usepackage{longtable}
\usepackage[russian]{babel}
\usepackage{color}

\theoremstyle{remark}

\renewcommand{\r}{\mathbb R}

\newcommand{\N}{\mathbb N}

\newcommand{\sign}{{\rm sign}}
\renewcommand{\le}{\leqslant}
\renewcommand{\ge}{\geqslant}

\renewcommand{\phi}{\varphi}

\renewcommand{\Re}{{\rm Re}\:}
\newcommand{\eqd}{\stackrel{d}{=}}
\newcommand{\pto}{\stackrel{P}{\longrightarrow}}

\title{{Convergence of random sums and statistics constructed from
samples with random sizes to the Linnik and Mittag-Leffler
distributions and their generalizations}\thanks{Research supported
by the Russian Science Foundation (project 14-11-00364).}}
\author{V. Yu. Korolev\thanks{Faculty of Computational Mathematics
and Cybernetics, Lomonosov Moscow State University; Institute of
Informatics Problems, Federal Research Center <<Computer Science and
Control>> of the Russian Academy of Sciences; vkorolev@cs.msu.ru},
A. I. Zeifman\thanks{Vologda State University; Institute of
Informatics Problems, Federal Research Center <<Computer Science and
Control>> of the Russian Academy of Sciences, ISEDT RAS;
a$\_$zeifman@mail.ru}}

\date{}

\begin{document}

\maketitle

{

\small

{\bf Abstract:} We present some product representations for random
variables with the Linnik, Mittag-Leffler and Weibull distributions
and establish the relationship between the mixing distributions in
these representations. Based on these representations, we prove some
limit theorems for a wide class of rather simple statistics
constructed from samples with random sized including, e. g., random
sums of independent random variables with finite variances, maximum
random sums, extreme order statistics, in which the Linnik and
Mittag-Leffler distributions play the role of limit laws. Thus we
demonstrate that the scheme of geometric summation is far not the
only asymptotic setting (even for sums of independent random
variables) in which the Mittag-Leffler and Linnik laws appear as
limit distributions. The two-sided Mittag-Leffler and one-sided
Linnik distribution are introduced and also proved to be limit laws
for some statistics constructed from samples with random sizes.

\smallskip


}





\section{Introduction}

Usually the Mittag-Leffler and Linnik distributions are mentioned in
the literature together as examples of geometric stable
distributions. Since these distributions are very often pointed at
as weak limits for geometric random sums, there might have emerged a
prejudice that the scheme of {\it geometric} summation is the only
asymptotic setting within which these distributions can be limiting
for sums of independent and identically distributed random
variables. This prejudice is accompanied by the suspicion that
non-trivial ($\delta<1$, $\alpha<2$) Mittag-Leffler and Linnik laws
can be limiting only for sums in which the summands have infinite
variances.

The aim of this paper is to dispel this prejudice by presenting some
examples of limit theorems for a wide class of rather simple
statistics constructed from samples with random sizes including, e.
g., random sums of independent random variables {\it with finite
variances}, maximum random sums, extreme order statistics in which
the Linnik and Mittag-Leffler distributions play the role of limit
laws. We will demonstrate that the scheme of geometric summation is
far not the only asymptotic setting (even for sums of independent
random variables!) in which the Mittag-Leffler and Linnik laws
appear as limit distributions.

The main tools used to prove the limit theorems in this paper are
mixture representations for the Linnik, Mittag-Leffler and Weibull
distributions also presented here. Some of these representations
were known (mixture representations for the Linnik and
Mittag-Leffler laws were the objects of investigation in
\cite{Devroye1990, ErdoganOstrovskii1997, ErdoganOstrovskii1998,
KotzOstrovskii1996, Pakes1992, Kozubowski1998, Kozubowski1999}),
some of them are new. These representations open the way to
establish the close analytic and asymptotic relations between these
two laws.

Another obvious reason for which the Mittag-Leffler and Linnik
distributions are often brought together in the literature is the
formal similarity of the Laplace transform of the former and the
Fourier--Stieltjes transform of the latter. We will show that
actually the link between these two laws (and some laws related to
them) is much more interesting than the formal coincidence of their
transforms. We develop the known results on mixture representability
of the Linnik and Mittag-Leffler distributions and prove some new
results of this kind, thus finding a tight and clear analytical link
between the Linnik, Mittag-Leffler, stable and related
distributions. For example, it turns out that the Linnik
distribution with parameter $\alpha$ is a scale mixture of the
normal distributions with the mixing distribution being the
Mittag-Leffler law with parameter $\delta=\alpha/2$. Product
representations for the random variables with the Linnik and
Mittag-Leffler distributions obtained in the previous works were
aimed at the construction of convenient algorithms for the computer
generation of pseudo-random variables with these distributions. At
the same time, mixture representation for the Linnik distribution as
a scale mixture of normals opens the way for the construction of a
random-sum central limit theorem with the Linnik distribution as the
limit law. Moreover, in the ``if and only if'' version of the
random-sum central limit theorem presented in this paper the
Mittag-Leffler distribution {\it must} be the limit law for the
normalized number of summands.

Strange as it may seem, the results concerning the possibility of
representation of the Linnik distribution as a scale mixture of
normals have never been explicitly presented in the literature in
full detail although the property of the Linnik distribution to be a
normal scale mixture is something almost obvious. Perhaps, the paper
\cite{KotzOstrovskii1996} is the closest to this conclusion and
exposes the representability of the Linnik law as a scale mixture of
Laplace distributions with the mixing distribution written out
explicitly.

Other examples of the results obtained here are the representations
of the Mittag-Leffler law as a scale mixture of Weibull or
half-normal distributions, based on which we prove theorems
establishing the conditions for the distributions of extreme order
statistics in samples with random sizes or maximum random sums of
independent random variables with finite variances to converge to
the Mittag-Leffler law.

The paper is organized as follows. Section 2 presents the
definitions and basic properties of the Linnik and Mittag-Leffler
distributions. Section 3 contains basic definitions and auxiliary
results. The proofs of our main results are purposely indirect and
essentially rely on some new mixture properties of the Weibull
distribution also presented in Section 3. In Section 4 we prove the
representability of the Linnik distribution as the scale mixture of
normal laws with the Mittag-Leffler mixing distribution. We use this
result to decribe the asymptotics of the tail behavior of the Linnik
distribution. Here we also obtain the representation of the Linnik
distribution as a scale mixture of the Laplace laws with the mixing
distribution explicitly determined as that of the ratio of two
independent random variables with the same strictly stable
distribution concentrated on the nonnegative halfline. We use this
representation together with the result of \cite{KotzOstrovskii1996}
to obtain a by-product corollary which is the explicit
representation of the distribution density of the ratio of two
independent positive strictly stable random variables, thus giving a
new proof of a result of \cite{Dovidio2010}. In Section 5 we prove
some representations of the Mittag-Leffler distribution as a mixed
exponential or a mixed half-normal law. In Section 6 we prove and
discuss some criteria (that is, necessary and sufficient conditions)
for the convergence of the distributions of rather simple statistics
constructed from samples with random sizes including, e. g., random
sums of independent random variables with finite variances, maximum
random sums, extreme order statistics, to the Linnik and
Mittag-Leffler laws. The asymptotic theory of extreme values in
samples with random sizes is well-developed. The basics of this
theory were presented, say, in \cite{Berman1964,
BarndorffNielsen1964, Mogyorodi1967, GnedenkoSenusiBereksi1982}. A
detailed review of this theory can be found in \cite{Galambos1984}.
Dealing with extreme order statistics in Section 6 we consider a
special but rather important case where the sample size is generated
by a doubly stochastic Poisson process and consider the asymptotic
behavior of the so-called max-compound Cox processes introduced and
studied in \cite{KorolevSokolov2008}. Here we also present two
examples of the construction of ``appropriate'' random indices
possessing the desired asymptotic properties. In Section 7 the
symmetric two-sided Mittag-Leffler distribution and the one-sided
Linnik distribution are introduced. Here we prove theorems stating
that these laws can also be limit distributions for statistics
constructed from samples with random sizes such as random sums of
independent random variables with finite variances, maximum random
sums or extreme order statistics.

\section{The Mittag-Leffler and Linnik distributions}

\subsection{The Mittag-Leffler distributions}

The Mittag-Leffler probability distribution is the distribution of a
nonnegative random variable $M_{\delta}$ whose Laplace transform is
$$
\psi_{\delta}(s)\equiv {\sf E}e^{-sM_{\delta}}=\frac{1}{1+\lambda
s^{\delta}},\ \ \ s\ge0,\eqno(1)
$$
where $\lambda>0$, $0<\delta\le1$. For simplicity, in what follows
we will consider the standard scale case and assume that
$\lambda=1$.

The origin of the term {\it Mittag-Leffler distribution} is due to
that the probability density corresponding to Laplace transform (1)
has the form
$$
f_{\delta}^{M}(x)=\frac{1}{x^{1-\delta}}\sum_{n=0}^{\infty}\frac{(-1)^nx^{\delta
n}}{\Gamma(\delta n+1)}=-\frac{d}{dx}E_{\delta}(-x^{\delta}),\ \ \
x\ge0,\eqno(2)
$$
where $E_{\delta}(z)$ is the Mittag-Leffler function with index
$\delta$ that is defined as the power series
$$
E_{\delta}(z)=\sum_{n=0}^{\infty}\frac{z^n}{\Gamma(\delta n+1)},\ \
\ \delta>0,\ z\in\mathbb{Z}.
$$
Here $\Gamma(s)$ is Euler's gamma-function,
$$
\Gamma(s)=\int_{0}^{\infty}z^{s-1}e^{-z}dz,\ \ \ s>0.
$$
The distribution function corresponding to density (2) will be
denoted $F_{\delta}^{M}(x)$.

With $\delta=1$, the Mittag-Leffler distribution turns into the
standard exponential distribution, that is, $F_1^{M}(x)=
[1-e^{-x}]\mathbf{1}(x\ge 0)$, $x\in\mathbb{R}$ (here and in what
follows the symbol $\mathbf{1}(C)$ denotes the indicator function of
a set $C$). But with $\delta<1$ the Mittag-Leffler distribution
density has the heavy power-type tail: from the well-known
asymptotic properties of the Mittag-Leffler function it can be
deduced that if $0<\delta<1$, then
$$
f_\delta^{M}(x)\sim \frac{\sin(\delta\pi)\Gamma(\delta+1)}{\pi
x^{\delta+1}}\eqno(3)
$$
as $x\to\infty$, see, e. g., \cite{Kilbas2014}.

It is well-known that the Mittag-Leffler distribution is stable with
respect to geometric summation (or {\it geometrically stable}). This
means that if $X_1,X_2,\ldots$ are independent random variables and
$V_p$ is the random variable independent of $X_1,X_2,\ldots$ and
having the geometric distribution
$$
{\sf P}(V_p=n)=p(1-p)^{n-1},\ \ \ n=1,2,\ldots,\ \ \
p\in(0,1),\eqno(4)
$$
then for each $p\in(0,1)$ there exists a constant $a_p>0$ such that
$a_p\big(X_1+\ldots+X_{V_p}\big)\Longrightarrow M_{\delta}$ as $p\to
0$, see, e. g., \cite{Bunge1996} or \cite{KlebanovRachev1996} (the
symbol $\Longrightarrow$ hereinafter denotes convergence in
distribution). Moreover, as far ago as in 1965 it was shown by
I.~Kovalenko \cite{Kovalenko1965} that the distributions with
Laplace transforms (1) are the only possible limit laws for the
distributions of appropriately normalized geometric sums of the form
$a_p\big(X_1+\ldots+X_{V_p}\big)$ as $p\to0$, where $X_1,X_2,\ldots$
are independent identically distributed {\it nonnegative} random
variables and $V_p$ is the random variable with geometric
distribution (4) independent of the sequence $X_1,X_2,\ldots$ for
each $p\in(0,1)$. The proofs of this result were reproduced in
\cite{GnedenkoKovalenko1968, GnedenkoKovalenko1989} and
\cite{GnedenkoKorolev1996}. In these books the class of
distributions with Laplace transforms (1) was not identified as the
class of Mittag-Leffler distributions but was called {\it class}
$\mathcal{K}$ after I.~Kovalenko.

Twenty five years later this limit property of the Mittag-Leffler
distributions was re-discovered by A.~Pillai in \cite{Pillai1989,
Pillai1990} who proposed the term {\it Mittag-Leffler distribution}
for the distribution with Laplace transform (1). Perhaps, since the
works \cite{Kovalenko1965, GnedenkoKovalenko1968,
GnedenkoKovalenko1989} were not easily available to probabilists,
the term {\it class $\mathcal{K}$ distribution} did not take roots
in the literature whereas the term {\it Mittag-Leffler distribution}
became conventional.

The Mittag-Leffler distributions are of serious theoretical interest
in the problems related to thinned (or rarefied) homogeneous flows
of events such as renewal processes or anomalous diffusion or
relaxation phenomena, see \cite{WeronKotulski1996,
GorenfloMainardi2006} and the references therein.

\subsection{The Linnik distributions}

In 1953 Yu. V. Linnik \cite{Linnik1953} introduced the class of
symmetric probability distributions defined by the characteristic
functions
$$
\mathfrak{f}^L_{\alpha}(t)=\frac{1}{1+|t|^{\alpha}},\ \ \
t\in\mathbb{R},\eqno(5)
$$
where $\alpha\in(0,2]$. Later the distributions of this class were
called {\it Linnik distributions} \cite{Kotz2001} or {\it
$\alpha$-Laplace distributions} \cite{Pillai1985}. In this paper we
will keep to the first term that has become conventional. With
$\alpha=2$, the Linnik distribution turns into the Laplace
distribution corresponding to the density
$$
f^{\Lambda}(x)=\textstyle{\frac12}e^{-|x|},\ \ \
x\in\mathbb{R}.\eqno(6)
$$
A random variable with Laplace density (6) and its distribution
function will be denoted $\Lambda$ and $F^{\Lambda}(x)$,
respectively.

The Linnik distributions possess many interesting analytic
properties such as unimodality \cite{Laha1961} and infinite
divisibility \cite{Devroye1990}, existence of an infinite peak of
the density for $\alpha\le1$ \cite{Devroye1990}, etc. In
\cite{KotzOstrovskiiHayfavi1995a, KotzOstrovskiiHayfavi1995b} a
detailed investigation of analytic and asymptotic properties of the
density of the Linnik distribution was carried out. However,
perhaps, most often Linnik distributions are recalled as examples of
geometric stable distributions.

A random variable with the Linnik distribution with parameter
$\alpha$ will be denoted $L_{\alpha}$. Its distribution function and
density will be denoted $F_{\alpha}^{L}$ and $f_{\alpha}^{L}$,
respectively. As this is so, from (5) and (6) it follows that
$F_2^{L}(x)\equiv F^{\Lambda}(x)$, $x\in\mathbb{R}$.

\section{Basic notation and auxiliary results}

Most results presented below actually concern special mixture
representations for probability distributions. However, without any
loss of generality, for the sake of visuality and compactness of
formulations and proofs we will represent the results in terms of
the corresponding random variables assuming that all the random
variables mentioned in what follows are defined on the same
probability space $(\Omega,\,\mathfrak{A},\,{\sf P})$.

The random variable with the standard normal distribution function
$\Phi(x)$ will be denoted $X$,
$$
{\sf
P}(X<x)=\Phi(x)=\frac{1}{\sqrt{2\pi}}\int_{-\infty}^{x}e^{-z^2/2}dz,\
\ \ \ x\in\mathbb{R}.
$$
Let $\Psi(x)$, $x\in\mathbb{R}$, be the distribution function of the
maximum of the standard Wiener process on the unit interval,
$\Psi(x)=2\Phi\big(\max\{0,x\}\big)-1$, $x\in\mathbb{R}$. It is easy
to see that $\Psi(x)={\sf P}(|X|<x)$. Therefore, sometimes $\Psi(x)$
is said to determine the {\it half-normal} distribution.

Throughout the paper the symbol $\eqd$ will denote the coincidence
of distributions.

The distribution function and the density of the strictly stable
distribution with the characteristic exponent $\alpha$ and shape
parameter $\theta$ defined by the characteristic function
$$
\mathfrak{g}_{\alpha,\theta}(t)=\exp\big\{-|t|^{\alpha}\exp\{-{\textstyle\frac12}i\pi\theta\alpha\,\mathrm{sign}t\}\big\},\
\ \ \ t\in\r,\eqno(7) 
$$
with $0<\alpha\le2$, $|\theta|\le\min\{1,\frac{2}{\alpha}-1\}$, will
be denoted by $G_{\alpha,\theta}(x)$ and $g_{\alpha,\theta}(x)$,
respectively (see, e. g., \cite{Zolotarev1983}). Any random variable
with the distribution function $G_{\alpha,\theta}(x)$ will be
denoted $S_{\alpha,\theta}$.

From (7) it follows that the characteristic function of a symmetric
($\theta=0$) strictly stable distribution has the form
$$
\mathfrak{g}_{\alpha,0}(t)=e^{-|t|^{\alpha}},\ \ \ t\in\r. \eqno(8) 
$$

\smallskip

{\sc Lemma 1}. {\it Let $\alpha\in(0,2]$, $\alpha'\in(0,1]$. Then}
$$
S_{\alpha\alpha',0}\eqd S_{\alpha,0}S_{\alpha',1}^{1/\alpha}
$$
{\it where the random variables on the right-hand side are
independent.}

\smallskip

{\sc Proof}. See, e. g., theorem 3.3.1 in \cite{Zolotarev1983}.

\smallskip

{\sc Corollary 1.} {\it A symmetric strictly stable distribution
with the characteristic exponent $\alpha$ is a scale mixture of
normal laws in which the mixing distribution is the one-sided
strictly stable law $(\theta=1)$ with the characteristic exponent
$\alpha/2:$
$$
S_{\alpha,0}\eqd X\sqrt{S_{\alpha/2,1}} \eqno(9)
$$
with the random variables on the right-hand side being independent.}

\smallskip

In terms of distribution functions the statement of corollary 1 can
be written as
$$
G_{\alpha,0}(x)=\int_{0}^{\infty}\Phi\Big(\frac{x}{\sqrt{z}}\Big)dG_{\alpha/2,1}(z),\
\ \ x\in\r.
$$

\smallskip

Let $\gamma>0$. The distribution of the random variable
$W_{\gamma}$:
$$
{\sf
P}\big(W_{\gamma}<x\big)=\big[1-e^{-x^{\gamma}}\big]\mathbf{1}(x\ge
0),\ \ \ x\in\mathbb{R},
$$
is called the {\it Weibull distribution} with shape parameter
$\gamma$. It is obvious that $W_1$ is the random variable with the
standard exponential distribution: ${\sf
P}(W_1<x)=\big[1-e^{-x}\big]{\bf 1}(x\ge0)$. The Weibull
distribution with $\gamma=2$, that is, ${\sf
P}(W_2<x)=\big[1-e^{-x^2}\big]{\bf 1}(x\ge0)$ is called the Rayleigh
distribution.

It is easy to see that if $\gamma>0$ and $\gamma'>0$, then ${\sf
P}(W_{\gamma'}^{1/\gamma}\ge x)={\sf P}(W_{\gamma'}\ge
x^{\gamma})=e^{-x^{\gamma\gamma'}}={\sf P}(W_{\gamma\gamma'}\ge x)$,
$x\ge 0$, that is, for any $\gamma>0$ and $\gamma'>0$
$$
W_{\gamma\gamma'}\eqd W_{\gamma'}^{1/\gamma}.\eqno(10) 
$$

It can be shown that each Weibull distribution with parameter
$\gamma\in(0,1]$ is a mixed exponential distribution. In order to
prove this we first make sure that each Weibull distribution with
parameter $\gamma\in(0,2]$ is a scale mixture of the Rayleigh
distributions.

For $\alpha\in(0,1]$ denote $T_{\alpha}=2S_{\alpha,1}^{-1}$, where
$S_{\alpha,1}$ is a random variable with one-sided strictly stable
density $g_{\alpha,1}(x)$.

\smallskip

{\sc Lemma 2.} {\it For any $\gamma\in(0,2]$ we have}
$$
W_{\gamma}\eqd W_2 \sqrt{T_{\gamma/2}},
$$
{\it where the random variables on the right-hand side are
independent.}

\smallskip

{\sc Proof}. Write relation (9) in terms of characteristic functions
with the account of (7):
$$
e^{-|t|^{\alpha}}=\int_{0}^{\infty}\exp\{-{\textstyle\frac12}t^2z\}g_{\alpha/2,1}(z)dz,
\ \ \ t\in\mathbb{R}.\eqno(11) %
$$
Formally letting $|t|=x$ in (11), where $x\ge0$ is an arbitrary
nonnegative number, we obtain
$$
{\sf
P}(W_{\gamma}>x)=e^{-x^{\gamma}}=\int_{0}^{\infty}\exp\{-{\textstyle\frac12}x^2z\}g_{\gamma/2,1}(z)dz.\eqno(12) 
$$
At the same time it is obvious that if $W_2$ and $S_{\gamma/2,1}$
are independent, then
$$
{\sf P}\big(W_2 \sqrt{T_{\gamma/2}}>x\big)={\sf
P}\big(W_2>x\sqrt{{\textstyle\frac12}S_{\gamma/2,1}}\big)=
\int_{0}^{\infty}\exp\{-{\textstyle\frac12}x^2z\}g_{\gamma/2,1}(z)dz.\eqno(13) 
$$
Since the right-hand sides of (12) and (13) coincide identically in
$x\ge0$, the left-hand sides of these relations coincide as well.
The lemma is proved.

\smallskip

{\sc Lemma 3}. {\it For any $\gamma\in(0,1]$, the Weibull
distribution with parameter $\gamma$ is a mixed exponential
distribution}:
$$
W_{\gamma}\eqd W_1 T_{\gamma}.\eqno(14) 
$$
{\it where the random variables on the right-hand side of $(14)$ are
independent.}

\smallskip

{\sc Proof.} It is easy to see that ${\sf P}(W_1^{1/\gamma}\ge
x)={\sf P}(W_1\ge x^{\gamma})=e^{-x^{\gamma}}={\sf P}(W_{\gamma}\ge
x)$, $x\ge 0$, that is,
$$
W_{\gamma}\eqd W_1^{1/\gamma}\eqno(15) 
$$
for any $\gamma>0$. From $(15)$ it follows that $W_2\eqd\sqrt{W_1}$.
Therefore, from lemma 2 it follows that for $\gamma\in(0,2]$ we have
$$
W_{\gamma}\eqd W_2 \sqrt{T_{\gamma/2}}\eqd\sqrt{W_1 T_{\gamma/2}}
$$
or, with the account of $(15)$,
$$
W_{\gamma/2}\eqd W_{\gamma}^2\eqd W_1 T_{\gamma/2}.
$$
Re-denoting $\gamma/2\longmapsto \gamma\in(0,1]$, we obtain the
desired assertion.

\smallskip

In \cite{Devroye1990} the following statement was proved. Here its
formulation is extended with the account of (10).

\smallskip

{\sc Lemma 4} \cite{Devroye1990}. {\it For any $\alpha\in(0,2]$, the
Linnik distribution with parameter $\alpha$ is a scale mixture of a
symmetric stable distribution with the Weibull mixing distribution
with parameter $\alpha/2$, that is,}
$$
L_{\alpha}\eqd S_{\alpha,0}W_{\alpha}\eqd
S_{\alpha,0}\sqrt{W_{\alpha/2}},
$$
{\it where the random variables on the right-hand side are
independent}.

\smallskip

{\sc Lemma 5.} {\it For any $\delta\in(0,1]$, the Mittag-Leffler
distribution with parameter $\delta$ is a scale mixture of a
one-sided stable distribution with the Weibull mixing distribution
with parameter $\delta/2$, that is,}
$$
M_{\delta}\eqd S_{\delta,1}W_{\delta}\eqd
S_{\delta,1}\sqrt{W_{\delta/2}},
$$
{\it where the random variables on the right-hand side are
independent}.

\smallskip

{\sc Proof}. This statement has already become folklore. For the
purpose of convenience we give its elementary proof without any
claims for priority. Let $S_{\delta,1}$ be a positive strictly
stable random variable. As is known, its Laplace transform is
$\psi(s)={\sf E}e^{-sS_{\delta,1}}=e^{-s^{\delta}}$, $s\ge0$. Then
with the account of (10) by the Fubini theorem the Laplace transform
of the product $S_{\delta,1}W_{\delta}$ is
$$
{\sf E}\exp\{-sS_{\delta,1}W_{\delta}\}= {\sf
E}\exp\{-sS_{\delta,1}W_1^{1/\delta}\}={\sf E}{\sf
E}\big(\exp\{-sS_{\delta,1}W_1^{1/\delta}\}\big|W_1\big)=
\int_{0}^{\infty}e^{-(sz^{1/\delta})^{\delta}}e^{-z}dz=
$$
$$
=\int_{0}^{\infty}e^{-z(s^{\delta}+1)}dz=\frac{1}{1+s^{\delta}}={\sf
E}e^{-sM_{\delta}},\ \ \ s\ge0.
$$
The lemma is proved.

\smallskip

Let $\rho\in(0,1)$. In \cite{Kozubowski1998} it was demonstrated
that the function
$$
f_{\rho}^{K}(x)=\frac{\sin(\pi\rho)}{\pi\rho[x^2+2x\cos(\pi\rho)+1]},\
\ \ x\in(0,\infty),\eqno(16)
$$
is a probability density on $(0,\infty)$. Let $K_{\rho}$ be a random
variable with density (16).

\smallskip

{\sc Lemma 6} \cite{Kozubowski1998}. {\it Let $0<\delta<\delta'\le1$
and $\rho=\delta/\delta'<1$. Then}
$$
M_{\delta}\eqd M_{\delta'}K_{\rho}^{1/\delta}
$$
{\it where the random variables on the right-hand side are
independent.}

\smallskip

With $\delta'=1$ we have

\smallskip

{\sc Corollary 2} \cite{Kozubowski1998}. {\it Let $0<\delta<1$. Then
the Mittag-Leffler distribution with parameter $\delta$ is mixed
exponential}:
$$
M_{\delta}\eqd K_{\delta}^{1/\delta}W_1
$$
{\it where the random variables on the right-hand side are
independent.}

\smallskip

Let $0<\alpha<\alpha'\le2$. In \cite{KotzOstrovskii1996} it was
shown that the function
$$
f_{\alpha,\alpha'}^{Q}(x)= \frac{\alpha'\sin(\pi\alpha/\alpha')
x^{\alpha-1}}{\pi[1+x^{2\alpha}+2x^{\alpha}\cos(\pi\alpha/\alpha')]},\
\ \ x>0,\eqno(17)
$$
is a probability density on $(0,\infty)$. Let $Q_{\alpha,\alpha'}$
be a random variable whose probability density is
$f_{\alpha,\alpha'}^{Q}(x)$.

\smallskip

{\sc Lemma 7} \cite{KotzOstrovskii1996}. {\it Let
$0<\alpha<\alpha'\le2$. Then}
$$
L_{\alpha}\eqd L_{\alpha'}Q_{\alpha,\alpha'},
$$
{\it where the random variables on the right-hand side are
independent}.

\smallskip

With $\alpha'=2$ we have

\smallskip

{\sc Corollary 3} \cite{KotzOstrovskii1996}. {\it Let $0<\alpha<2$.
Then the Linnik distribution with parameter $\alpha$ is a scale
mixture of Laplace distributions corresponding to density $(5)$}:
$$
L_{\alpha}\eqd \Lambda Q_{\alpha,2}
$$
{\it where the random variables on the right-hand side are
independent.}

\smallskip

For the sake of completeness, we will demonstrate that the Weibull
distributions possess the same property as the Linnik and
Mittag-Leffler distributions presented in lemmas 6 and 7: any
distribution of the corresponding class can be represented as a
scale mixture of a distribution from the same class with larger
parameter.

Relation (14) implies the following statement generalizing lemmas 2
and 3 and stating that the Weibull distribution with an arbitrary
positive shape parameter $\gamma$ is a scale mixture of the Weibull
distribution with an arbitrary positive shape parameter
$\gamma'>\gamma$.

\smallskip

{\sc Lemma 8}. {\it Let $\gamma'>\gamma>0$ be arbitrary numbers.
Then
$$
W_{\gamma}\eqd W_{\gamma'}\cdot T_{\alpha}^{1/\gamma'},
$$
where $\alpha=\gamma/\gamma'\in(0,1)$ and the random variables on
the right-hand side are independent}.

\smallskip

{\sc Proof}. In lemma 3 we showed that a Weibull distribution with
parameter $\alpha\in(0,1]$ is a mixed exponential distribution.
Indeed, from (14) it follows that
$$
e^{-x^{\alpha}}={\sf P}(W_{\alpha}>x)={\sf
P}(W_1>{\textstyle\frac12}S_{\alpha,1}x)=
\int_{0}^{\infty}e^{-\frac12 zx}g_{\alpha,1}(z)dz,\ \ \ x\ge0.
$$
Therefore, for any $\gamma'>\gamma>0$, denoting
$\alpha=\gamma/\gamma'$ (as this is so, $\alpha\in(0,1)$), for any
$x\in\mathbb{R}$ we obtain
$$
{\sf P}(W_{\gamma}>x)=e^{-x^{\gamma}}=e^{-x^{\gamma'\alpha}}={\sf
P}(W_{\alpha}>x^{\gamma'})={\sf
P}(W_1>{\textstyle\frac12}S_{\alpha,1}x^{\gamma'})=
$$
$$
=\int_{0}^{\infty}e^{-\frac12
zx^{\gamma'}}g_{\alpha,1}(z)dz=\int_{0}^{\infty}{\sf
P}\big(W_{\gamma'}>x({\textstyle\frac{1}{2}}z)^{1/\gamma'}\big)g_{\alpha,1}(z)dz={\sf
P}(W_{\gamma'}\cdot T_{\alpha}^{1/\gamma'}>x),
$$
The lemma is proved.

\smallskip

It should be noted that if $0<\gamma<\gamma'<2$, then the assertion
of lemma 8 directly follows from theorem 3.3.1 of
\cite{Zolotarev1983} due to the formal coincidence of the
characteristic function of a strictly stable law and the
complementary Weibull distribution function (see the proof of lemma
2).

\smallskip

{\sc Corollary 4.} {\it Let $\gamma\ge1$ be arbitrary. Then the
exponential distribution is a scale mixture of the Weibull laws with
parameter $\gamma:$
$$
W_1\eqd W_{\gamma}\cdot T_{1/\gamma}^{1/\gamma},
$$
where the random variables on the right-hand side are independent.}

\section{Representation of the Linnik distribution as a scale
mixture of normal or Laplace distributions and related results}

\subsection{The representation of the Linnik distribution as a normal scale mixture}

In all the products of random variables mentioned below the
multipliers are assumed independent.

\smallskip

{\sc Theorem 1.} {\it Let $\alpha\in(0,2]$, $\alpha'\in(0,1]$. Then}
$$
L_{\alpha\alpha'}\eqd S_{\alpha,0}M_{\alpha'}^{1/\alpha}.
$$

\smallskip

{\sc Proof}. From lemma 4 we have
$$
L_{\alpha\alpha'}\eqd
S_{\alpha\alpha',0}\sqrt{W_{\alpha\alpha'/2}}.\eqno(18)
$$
Continuing (18) with the account of lemma 1, we obtain
$$
L_{\alpha\alpha'}\eqd
S_{\alpha,0}S_{\alpha',1}^{1/\alpha}\sqrt{W_{\alpha\alpha'/2}}.\eqno(19)
$$
From (10) and lemma 5 it follows that
$$
S_{\alpha',1}^{1/\alpha}\sqrt{W_{\alpha\alpha'/2}}\eqd
S_{\alpha',1}^{1/\alpha}W_{\alpha\alpha'}\eqd
(S_{\alpha',1}W_{\alpha'})^{1/\alpha}\eqd M_{\alpha'}^{1/\alpha}.
$$
The theorem is proved.

\smallskip

As far as we know, the following result has never been explicitly
presented in the literature in full detail although the property of
the Linnik distribution to be a normal scale mixture is something
almost obvious.

\smallskip

{\sc Corollary 5}. {\it For each $\alpha\in(0,2]$, the Linnik
distribution with parameter $\alpha$ is the scale mixture of
zero-mean normal laws with mixing Mittag-Leffler distribution with
twice less parameter $\alpha/2$}:
$$
L_{\alpha}\eqd X\sqrt{M_{\alpha/2}},\eqno(20)
$$
{\it where the random variables on the right-hand side are
independent}.

\subsection{The tail behavior of the Linnik distribution}

From $(20)$ we can easily characterize the tail behavior of the
Linnik distribution. For this purpose we will use the following
statement proved in \cite{AntonovKoksharov2006}.

\smallskip

{\sc Lemma 9} \cite{AntonovKoksharov2006}. {\it If a distribution
function $F(x)$ has the form
$$
F(x)=\int_{0}^{\infty}\Phi\Big(\frac{x}{\sqrt{u}}\Big)dG(u),\ \ \
x\in\mathbb{R},
$$
where $G(u)$ is a distribution function such that $G(0)=0$, and
$\rho$ and $C$ are positive numbers, then the conditions
$$
\limsup_{x\to\infty}x^{\rho}[1-F(x)]=C
$$
and
$$
\limsup_{u\to\infty}u^{\rho/2}[1-G(u)]=2C
$$
are equivalent.}

\smallskip

From lemma 9 with the account of $(3)$ and $(20)$ we obtain the
following statement.

\smallskip

{\sc Corollary 6.} {\it The tail behavior of the Linnik distribution
$L_{\alpha}$ with parameter $\alpha\in(0,2)$ as $x\to\infty$ is
described by the relation}
$$
\limsup_{x\to\infty}x^{\alpha/2}[1-L_{\alpha}(x)]=\frac{\alpha}{2\pi}\sin\Big(\frac{\alpha\pi}{2}\Big)\Gamma\Big(\frac{\alpha}{2}+1\Big).
$$

\smallskip

In other words, the above reasoning gives one more proof that if
$0<\alpha<2$, then $1-L_{\alpha}(x)=O(x^{-\alpha/2})$ as
$x\to\infty$, not involving the tail properties of stable
distributions.

\subsection{The representation of the Linnik distribution as a scale mixture of Laplace distributions}

It should be noted that by lemma 5 representation (20) can be
rewritten as
$$
L_{\alpha}\eqd X\sqrt{S_{\alpha/2,1}W_{\alpha/2}}.\eqno(21)
$$
From lemma 3 we have
$$
W_{\alpha/2}\eqd \frac{2W_1}{S'_{\alpha/2,1}}.
$$
Hence, from (21) it follows that
$$
L_{\alpha}\eqd X\sqrt{2W_1\frac{S_{\alpha/2,1}}{S'_{\alpha/2,1}}}
$$
where the independent random variables $S_{\alpha/2,1}$ and
$S'_{\alpha/2,1}$ have one and the same one-sided strictly stable
distribution with characteristic exponent $\alpha/2$ and are
independent of the exponentially distributed random variable $W_1$.
It is well known that
$$
X\sqrt{2W_1}\eqd \Lambda\eqno(22)
$$
(see, e. g., the example on p. 272 of \cite{BeningKorolev2002} or
lemma 10 below). Therefore we obtain one more mixture representation
for the Linnik distribution.

\smallskip

{\sc Theorem 2}. {\it For each $\alpha\in(0,2]$, the Linnik
distribution with parameter $\alpha$ is the scale mixture of the
Laplace laws corresponding to density $(6)$ with mixing distribution
being that of the ratio of two independent random variables having
one and the same one-sided strictly stable distribution with
characteristic exponent $\alpha/2$}:
$$
L_{\alpha}\eqd \Lambda\sqrt{\frac{S_{\alpha/2,1}}{S'_{\alpha/2,1}}},
$$
{\it where the random variables on the right-hand side are
independent}.

\smallskip

It is easy to see that scale mixtures of Laplace distribution (6)
are identifiable, that is, if
$$
\Lambda Y\eqd\Lambda Y'
$$
where $Y$ and $Y$ are nonnegative random variables independent of
$\Lambda$, then $Y\eqd Y'$. Indeed, with the account of (22), the
last relation turns into
$$
X\sqrt{2W_1 Y^2}\eqd X\sqrt{2W_1 (Y')^2},\eqno(23)
$$
where the random the multipliers on both sides are independent. But,
as is known, scale mixtures of zero-mean normals are identifiable
(see \cite{Teicher1961}). Therefore, (23) implies that
$$
W_1 Y^2\eqd W_1(Y')^2.\eqno(24)
$$
The complementary mixed exponential distribution functions of the
random variables related by (24) are the Laplace transforms of $Y^2$
and $(Y')^2$, respectively. Relation (24) means that these Laplace
transforms identically coincide:
$$
\int_{0}^{\infty}e^{-sz}d{\sf
P}(Y^2<z)\equiv\int_{0}^{\infty}e^{-sz}d{\sf P}\big((Y')^2<z\big),\
\ \ s\ge0.
$$
Hence, the distributions of the random variables $Y^2$ and $(Y')^2$
coincide and hence, the distributions of $Y$ and $Y'$ coincide as
well since these random variables were originally assumed
nonnegative.

\subsection{Some properties of the mixing distributions}

Comparing the statement of theorem 2 with the assertion of corollary
3 with the account of identifiability of scale mixtures of Laplace
distributions (6) we arrive at the relation
$$
Q_{\alpha,2}\eqd\sqrt{\frac{S_{\alpha/2,1}}{S'_{\alpha/2,1}}}.\eqno(25)
$$
The combination of (17) and (25) gives one more, possibly simpler,
proof of the following by-product result concerning the properties
of stable distributions obtained in \cite{Dovidio2010}. This result
offers an explicit representation for the density of the ratio of
two independent stable random variables in terms of elementary
functions although with the exception of one case, the L{\'e}vy
distribution $(\alpha=\frac12)$, such representations for the
densities of nonnegative stable random variables themselves do not
exist.

\smallskip

{\sc Corollary 7}. {\it Let $S_{\alpha,1}$ and $S'_{\alpha,1}$ be
two independent random variables having one and the same one-sided
strictly stable distribution with characteristic exponent
$\alpha\in(0,1)$. Then the probability density $p_{\alpha}(x)$ of
the ratio $S_{\alpha,1}/S'_{\alpha,1}$ has the form
$p_{\alpha}(x)=2xf_{2\alpha,2}^{Q}(x^2),$ $x\ge0$, where
$f_{2\alpha,2}^{Q}$ was defined in $(17)$, that is,}
$$
p_{\alpha}(x)=\frac{4\sin(\pi\alpha)
x^{4\alpha-1}}{\pi[1+x^{8\alpha}+2x^{4\alpha}\cos(\pi\alpha)]},\ \ \
x>0.\eqno(26)
$$

\smallskip

By the reasoning similar to that used to prove corollary 6 we can
obtain the following relation linking the distributions of the
random variables $K_{\delta}$ and $Q_{\alpha,2}$: for any
$\delta\in(0,1)$
$$
K_{\delta}^{1/\delta}\eqd
2\cdot\frac{S_{\delta,1}}{S'_{\delta,1}}\eqd
2Q^2_{2\delta,2}.\eqno(27)
$$
Hence, $K_{\delta}\eqd 2^{\delta}Q^{2\delta}_{2\delta,2}$.

Using lemmas 1, 4 and 5 it  it is possible to obtain more product
representations for the Mittag-Leffler- and Linnik-distributed
random variables and hence, more mixture representations for these
distributions.

\section{Exponential and half-normal mixture representations for the Mittag-Leffler distribution}

\subsection{The Mittag-Leffler distribution as a mixed exponential distribution}

As concerns the Mittag-Leffler distribution, from lemmas 3 and 5 we
obtain the following statement analogous to theorem 2.

\smallskip

{\sc Theorem 3}. {\it For each $\delta\in(0,1]$, the Mittag-Leffler
distribution with parameter $\delta$ is the mixed exponential
distribution with mixing distribution being that of twice the ratio
of two independent random variables having one and the same
one-sided strictly stable distribution with characteristic exponent
$\delta$}:
$$
M_\delta\eqd 2W_1\frac{S_{\delta,1}}{S'_{\delta,1}},
$$
{\it where the random variables on the right-hand side are
independent}.

\smallskip

From theorem 3 and corollary 7 we obtain the following
representation of the Mittag-Leffler distribution function
$F_{\delta}^{M}(x)$:
$$
F_{\delta}^{M}(x)=1-\frac{4\sin(\pi\delta)}{\pi}\int_{0}^{\infty}\frac{
z^{4\delta-1}e^{-2zx}dz}{1+z^{8\delta}+2z^{4\delta}\cos(\pi\delta)},\
\ \ x>0,\eqno(28)
$$
whence for the Mittag-Leffler density $f_{\delta}^{M}(x)$ we obtain
the integral representation
$$
f_{\delta}^{M}(x)=\frac{8\sin(\pi\delta)}{\pi}\int_{0}^{\infty}\frac{
z^{4\delta}e^{-2zx}dz}{1+z^{8\delta}+2z^{4\delta}\cos(\pi\delta)},\
\ \ x>0.
$$
A representation for the Linnik distribution similar to (28) was
obtained in \cite{Kozubowski1998}. We will use that representation
in Section 7.4.

\subsection{The Mittag-Leffler distribution as a mixture of half-normal distributions}

For a more thorough analysis of properties of the Mittag-Leffler law
as the limit distribution for random sums of independent random
variables we need the product representation of an exponential
random variable presented in what follows.

\smallskip

{\sc Lemma 10.} {\it The exponential distribution is a scale mixture
of half-normal laws. Namely, the relation
$$
W_1\eqd\sqrt{2W_1} |X|\eqno(29)
$$
holds, where the random variables on the right-hand side are
independent.}

\smallskip

{\sc Proof.} For $x>0$ we have
$$
{\sf P}\big(|X|\sqrt{W_1}<x\big)={\sf
E}\Psi\big(x/\sqrt{W_1}\big)=2{\sf E}\Phi\big(x/\sqrt{W_1}\big)-1=
2\int_{0}^{\infty}\Phi\big(x/\sqrt{z}\big)d[1-e^{-z}]-1=
$$
$$
=2\int_{0}^{\infty}\Big[\frac{1}{2}+
\frac{1}{\sqrt{2\pi}}\int_{0}^{x/\sqrt{z}}e^{-u^2/2}du\Big]e^{-z}dz-1=
\frac{\sqrt{2}}{\sqrt{\pi}}\int_{0}^{\infty}\int_{0}^{x/\sqrt{z}}e^{-u^2/2-z}dudz=
$$
$$=
\frac{\sqrt{2}}{\sqrt{\pi}}\int_{0}^{\infty}\int_{0}^{x^2/u^2}e^{-z}dz
e^{-u^2/2}du=
\frac{\sqrt{2}}{\sqrt{\pi}}\int_{0}^{\infty}\Bigl(
1-\exp\Bigl\{-\frac{x^2}{u^2}\Bigr\}\Bigr)e^{-u^2/2}du=
$$
$$=
1-\frac{\sqrt{2}}{\sqrt{\pi}}\int_{0}^{\infty}
\exp\Bigl\{-\frac{u^2}{2}-\frac{x^2}{u^2}\Bigr\}du=
1-e^{-\sqrt{2}x}={\sf P}(W_1<\sqrt{2}x),
$$
see, e. g., \cite{GradsteinRyzhik1971}, formula 3.325. This is
nothing else than the exponential distribution with parameter
$\sqrt{2}$. The lemma is proved.

\smallskip

From theorem 3 and lemma 10 we obtain the following representation
of the Mittag-Leffler distribution as a scale mixture of half-normal
laws.

\smallskip

{\sc Theorem 4.} {\it For $\delta\in(0,1]$ the Mittag-Leffler
distribution with parameter $\delta$ is a scale mixture of
half-normal laws}:
$$
M_{\delta}\eqd
|X|\sqrt{8W_1\Big(\frac{S_{\delta,1}}{S'_{\delta,1}}\Big)^2}.\eqno(30)
$$

\smallskip

In the subsequent sections an important role will be played by the
distribution which is mixing in (30). Denote
$$
H_{\delta}(x)={\sf
P}\bigg(W_1\Big(\frac{S_{\delta,1}}{S'_{\delta,1}}\Big)^2<\frac{x}{8}\bigg),\
\ \ x\ge0,\eqno(31)
$$
so that the assertion of theorem 4 can be written as
$$
F_{\delta}^{M}(x)=\int_{0}^{\infty}\Psi\Big(\frac{x}{\sqrt{u}}\Big)dH_{\delta}(u)=
2\int_{0}^{\infty}\Phi\Big(\frac{x}{\sqrt{u}}\Big)dH_{\delta}(u)-1,\
\ \ x\ge0.\eqno(32)
$$

With the account of corollary 7, for $0<\delta<1$ the density
corresponding to the distribution function $H_{\delta}(x)$ can be
written as
$$
h_{\delta}(x)=\frac{d}{dx}H_{\delta}(x)=\frac{d}{dx}{\sf
P}\bigg(W_1<\frac{x}{8}\Big(\frac{S'_{\delta,1}}{S_{\delta,1}}\Big)^2\bigg)=
\frac{\sin(\pi\delta)}{4\pi}\int_{0}^{\infty}\frac{z^{2\delta-1/2}e^{-xz/8}dz}{1+z^{4\delta}+2z^{2\delta}\cos(\pi\delta)},
\ \ \ x\ge0.
$$
If $\delta=1$, then obviously $H_1(x)=1-e^{-x/8}$, $x\ge0$.

From (32) and (3) it follows that if $0<\delta<1$, then, as
$x\to\infty$,
$$
1-F_{\delta}^{M}(x)\sim
\frac{\delta\sin(\delta\pi)\Gamma(\delta+1)}{\pi x^{\delta}}.
$$
Therefore, by lemma 9 from theorem 4 we obtain that
$$
\limsup_{x\to\infty}x^{\delta/2}[1-H_{\delta}(x)]=\frac{\delta\sin(\delta\pi)\Gamma(\delta+1)}{\pi},
$$
that is, if $0<\delta<1$, then $1-H_{\delta}(x)=O(x^{-\delta/2})$ as
$x\to\infty$.

\section{Convergence of the distributions of random sums and statistics
constructed from samples with random sizes to the Linnik and
Mittag-Leffler distributions}

\subsection{Convergence of the distributions of random sums to the Linnik
distribution}

In applied probability it is a convention that a model distribution
can be regarded as well-justified or adequate, if it is an {\it
asymptotic approximation}, that is, if there exists a rather simple
limit setting (say, schemes of maximum or summation of random
variables) and the corresponding limit theorem in which the model
under consideration manifests itself as a limit distribution. The
existence of such limit setting can provide a better understanding
of real mechanisms that generate observed statistical regularities.

As it has already been noted in the introduction, both the
Mittag-Leffler and Linnik laws are geometrically stable and are
therefore limit distributions for {\it geometric} random sums. In
this and subsequent sections we will demonstrate that the scheme of
geometric summation is far not the only asymptotic setting (even for
sums of independent random variables!) in which the Mittag-Leffler
and Linnik laws appear as limit distributions.

Product representations for the random variables with the Linnik and
Mittag-Leffler distributions obtained in the previous works were
aimed at the construction of convenient algorithms for the computer
generation of pseudo-random variables with these distributions. The
mixture representation for the Linnik distribution as a scale
mixture of normals obtained in corollary 4 opens the way for the
construction in this section of a random-sum central limit theorem
with the Linnik distribution as the limit law. Moreover, in this
``if and only if'' version of the random-sum central limit theorem
the Mittag-Leffler distribution {\it must} be the limit law for the
normalized number of summands.

Recall that the symbol $\Longrightarrow$ denotes the convergence in
distribution.

Consider independent not necessarily identically distributed random
variables $X_1,X_2,\ldots $ with ${\sf E}X_i=0$ and
$0<\sigma^2_i={\sf D}X_i<\infty$, $i\ge1$. For $n\in\N$ denote
$$
S^*_n=X_1+\ldots +X_n,\ \ \ \  B^2_n=\sigma^2_1+\ldots +\sigma^2_n.
$$
Assume that the random variables $X_1,X_2,\ldots $ satisfy the
Lindeberg condition: for any $\tau>0$
$$
\lim_{n\to\infty}\frac{1}{B^2_n}\sum_{i=1}^{n}\int_{|x|\ge\tau
B_n}^{} x^2d{\sf P}(X_i<x)=0.\eqno(33)
$$
It is well known that under these assumptions
$$
{\sf P}\big(S^*_n<B_nx\big)\Longrightarrow \Phi(x)
$$
(this is the classical Lindeberg central limit theorem).

Let $N_1,N_2,\ldots$ be a sequence of integer-valued nonnegative
random variables defined on the same probability space so that for
each $n\in\N$ the random variable $N_n$ is independent of the
sequence $X_1,X_2,\ldots$ Denote $S^*_{N_n}=X_1+\ldots +X_{N_n}$.
For definiteness, in what follows we assume that $\sum_{j=1}^0=0$.

Recall that a random sequence $N_1,N_2,\ldots$ is said to infinitely
increase in probability ($N_n\pto\infty$), if ${\sf P}(N_n\le
m)\longrightarrow 0$ as $n\to\infty$ for any $m\in(0,\infty)$.

Let $\{d_n\}_{n\ge1}$ be an infinitely increasing sequence of
positive numbers.

The proof of the main result of this section is based on the
following version of the random-sum central limit theorem.

\smallskip

{\sc Lemma 11} \cite{Korolev1994}. {\it Assume that the random
variables $X_1,X_2,\ldots$ and $N_1,N_2,\ldots$ satisfy the
conditions specified above. In particular, let Lindeberg condition
$(33)$ hold. Moreover, let $N_n\pto\infty$ as $n\to\infty$. A
distribution function $F(x)$ such that
$$
{\sf P}\Big(\frac{S^*_{N_n}}{d_n}<x\Big) \Longrightarrow F(x)
$$
as $n\to\infty$ exists if and only if there exists a distribution
function $H(x)$ satisfying the conditions
$$
H(0)=0,\ \ \
F(x)=\int_{0}^{\infty}\Phi\Big(\frac{x}{\sqrt{y}}\Big)dH(y),\ \
x\in\mathbb{R},
$$
and ${\sf P}(B^2_{N_n}<xd_n^2)\Longrightarrow H(x)$ $(n\to\infty)$.
}

\smallskip

{\sc Proof}. This statement is a particular case of a result proved
in \cite{Korolev1994}, also see theorem 3.3.2 in
\cite{GnedenkoKorolev1996}.

\smallskip

The following theorem gives a criterion (that is, {\it necessary and
sufficient} conditions) of the convergence of the distributions of
random sums of independent identically distributed random variables
with {\it finite} variances to the Linnik distribution.

\smallskip

{\sc Theorem 4.} {\it Let $\alpha\in(0,2]$. Assume that the random
variables $X_1,X_2,\ldots$ and $N_1,N_2,\ldots$ satisfy the
conditions specified above. In particular, let Lindeberg condition
$(33)$ hold. Moreover, let $N_n\pto\infty$ as $n\to\infty$. Then the
distributions of the normalized random sums $S^*_{N_n}$ converge to
the Linnik law with parameter $\alpha$, that is,
$$
{\sf P}\Big(\frac{S^*_{N_n}}{d_n}<x\Big) \Longrightarrow
F_{\alpha}^{L}(x)
$$
with some $d_n>0$, $d_n\to\infty$ as $n\to\infty$, if and only if
$$
\frac{B^2_{N_n}}{d^2_n}\Longrightarrow M_{\alpha/2}\ \ \
(n\to\infty).
$$
}

\smallskip

{\sc Proof}. This statement is a direct consequence of corollary 4
and lemma 11 with $H(x)=F_{\alpha/2}^{M}(x)$.

\smallskip

Note that if the random variables $X_1,X_2,\ldots$ are identically
distributed, then $\sigma_i=\sigma$, $i\in\N$, and the Lindeberg
condition holds automatically. In this case it is reasonable to take
$d_n=\sigma\sqrt{n}$. Hence, from theorem 4 in this case it follows
that for the convergence
$$
\frac{S^*_{N_n}}{\sigma\sqrt{n}}\Longrightarrow L_{\alpha}
$$
to hold as $n\to\infty$ it is necessary and sufficient that
$$
\frac{N_n}{n}\Longrightarrow M_{\alpha/2}.
$$

One more remark is that with $\alpha=2$ Theorem 4 involves the case
of convergence to the Laplace distribution.

\subsection{Convergence of the distributions of statistics
constructed from samples with random sizes to the Linnik
distribution}

In classical problems of mathematical statistics, the size of the
available sample, i. e., the number of available observations, is
traditionally assumed to be deterministic. In the asymptotic
settings it plays the role of infinitely increasing {\it known}
parameter. At the same time, in practice very often the data to be
analyzed is collected or registered during a certain period of time
and the flow of informative events each of which brings a next
observation forms a random point process. Therefore, the number of
available observations is unknown till the end of the process of
their registration and also must be treated as a (random)
observation. For example, this is so in insurance statistics where
during different accounting periods different numbers of insurance
events (insurance claims and/or insurance contracts) occur and in
high-frequency financial statistics where the number of events in a
limit order book during a time unit essentially depends on the
intensity of order flows. Moreover, contemporary statistical
procedures of insurance and financial mathematics do take this
circumstance into consideration as one of possible ways of dealing
with heavy tails. However, in other fields such as medical
statistics or quality control this approach has not become
conventional yet although the number of patients with a certain
disease varies from month to month due to seasonal factors or from
year to year due to some epidemic reasons and the number of failed
items varies from lot to lot. In these cases the number of available
observations as well as the observations themselves are unknown
beforehand and should be treated as random to avoid underestimation
of risks or error probabilities.

Therefore it is quite reasonable to study the asymptotic behavior of
general statistics constructed from samples with random sizes for
the purpose of construction of suitable and reasonable asymptotic
approximations. As this is so, to obtain non-trivial asymptotic
distributions in limit theorems of probability theory and
mathematical statistics, an appropriate centering and normalization
of random variables and vectors under consideration must be used. It
should be especially noted that to obtain reasonable approximation
to the distribution of the basic statistics, both centering and
normalizing values should be non-random. Otherwise the approximate
distribution becomes random itself and, for example, the problem of
evaluation of quantiles or significance levels becomes senseless.

In asymptotic settings, statistics constructed from samples with
random sizes are special cases of random sequences with random
indices. The randomness of indices usually leads to that the limit
distributions for the corresponding random sequences are
heavy-tailed even in the situations where the distributions of
non-randomly indexed random sequences are asymptotically normal see,
e. g., \cite{BeningKorolev2002, BeningKorolev2005,
GnedenkoKorolev1996}. For example, if a statistic which is
asymptotically normal in the traditional sense, is constructed on
the basis of a sample with random size having negative binomial
distribution, then instead of the expected normal law, the Student
distribution with power-type decreasing heavy tails appears as an
asymptotic law for this statistic.

Consider a problem setting that is traditional for mathematical
statistics. Let random variables $N_1,N_2,\ldots,X_1,X_2,\ldots,$ be
defined on one and the same probability space $(\Omega,{\cal A},
{\sf P})$. Assume that for each $n\ge 1$ the random variable $N_n$
takes only natural values and is independent of the sequence
$X_1,X_2,\ldots$ Let $T_n=T_n(X_1,\ldots,X_n)$ be a statistic, that
is, a measurable function of $X_1,\ldots,X_n$. For every $n\ge1$
define the random variable $T_{N_n}$ as
$$
T_{N_n}(\omega)=
T_{N_n(\omega)}\left(X_1(\omega),\ldots,X_{N_n(\omega)}(\omega)\right)
$$
for each $\omega\in\Omega$. As usual, the symbol $\Longrightarrow$
denotes convergence in distribution.

A statistic $T_n$ is said to be {\it asymptotically normal}, if
there exist $\delta>0$ and $\theta\in\r$ such that
$$
{\sf
P}\left(\delta\sqrt{n}\bigl(T_n-\theta\bigr)<x\right)\Longrightarrow\Phi(x)
\ \ \ (n\to\infty).\eqno(34)
$$

\smallskip

{\sc Lemma 12} \cite{Korolev1995}. {\it Assume that
$N_n\longrightarrow\infty$ in probability. Let the statistic $T_n$
be asymptotically normal in the sense of $(34)$. A distribution
function $F(x)$ such that
$$
{\sf P}\left(\delta\sqrt{n}\bigl(T_{N_n}-\theta\bigr)<x\right)
\Longrightarrow F(x)\ \ \ (n\to\infty),
$$
exists if and only if there exists a distribution function $H(x)$
satisfying the conditions
$$
H(0)=0,\ \ \ F(x)=\int_{0}^{\infty}\Phi\big(x\sqrt{y}\big)dH(y),\ \
x\in\mathbb{R},\ \ \ {\sf P}(N_n<nx)\Longrightarrow H(x) \
(n\to\infty).
$$
}

\smallskip

The following theorem gives a criterion (that is, {\it necessary and
sufficient} conditions) of the convergence of the distributions of
statistics, which are suggested to be asymptotically normal in the
traditional sense but are constructed from samples with random
sizes, to the Linnik distribution.

\smallskip

{\sc Theorem 5.} {\it Let $\alpha\in(0,2]$. Assume that the random
variables $X_1,X_2,\ldots$ and $N_1,N_2,\ldots$ satisfy the
conditions specified above and, moreover, let $N_n\pto\infty$ as
$n\to\infty$. Let the statistic $T_n$ be asymptotically normal in
the sense of $(34)$. Then the distribution of the statistic
$T_{N_n}$ constructed from samples with random sizes $N_n$ converges
to the Linnik law $F^L_{\alpha}(x)$ as $n\to\infty$, that is,
$$
{\sf P}\left(\delta\sqrt{n}\bigl(T_{N_n}-\theta\bigr)<x\right)
\Longrightarrow F^L_{\alpha}(x),
$$
if and only if
$$
\frac{N_n}{n}\Longrightarrow M_{\alpha/2}^{-1}\ \ \
(n\to\infty).\eqno(35)
$$
}

\smallskip

{\sc Proof}. This statement is a direct consequence of corollary 5
and lemma 12 with $H(x)={\sf P}(M_{\alpha/2}^{-1}<x)$.

\smallskip

From (28) and the absolute continuity of the Mittag-Leffler
distribution it follows that condition (35) can be written as
$$
\sup_{x>0}\bigg|{\sf
P}(N_n<nx)-\frac{4\sin(\pi\alpha/2)}{\pi}\int_{0}^{\infty}\frac{
z^{2\alpha-1}e^{-2z/x}dz}{1+z^{4\alpha}+2z^{2\alpha}\cos(\pi\alpha/2)}\bigg|=0.
$$

\subsection{Convergence of the distributions of extreme order
statistics constructed from samples with random sizes to the
Mittag-Leffler distribution}

Using lemmas 5 and 6 we can obtain the following representation of
the Mittag-Leffler distribution as a mixed Weibull distribution: if
$0<\delta<\delta'\le1$, then
$$
M_{\delta}\eqd
W_{\delta'}S_{\delta',1}K_{\delta/\delta'}^{1/\delta}.\eqno(36)
$$
It is well known that the Weibull distribution is a limit law for
extreme order statistics under an appropriate linear normalization.
This fact together with (36) open the way to prove that the
Mittag-Leffler distribution can be limiting for extreme order
statistics constructed from samples with random sizes.

In the book \cite{GnedenkoKorolev1996} it was proposed to model the
evolution of non-homogeneous chaotic stochastic processes, in
particular, the dynamics of financial markets by compound doubly
stochastic Poisson processes (compound Cox processes). This approach
got further grounds and development, say, in
\cite{BeningKorolev2002, Zeifman2015}. According to this approach
the flow of informative events, each of which generates the next
observation, is described by the stochastic point process $P(U(t))$
where $P(t)$, $t\geq0$, is a homogeneous Poisson process with unit
intensity and $U(t)$, $t\geq0$, is a random process independent of
$P(t)$ and possessing the properties: $U(0)=0$, ${\sf
P}(U(t)<\infty)=1$ for any $t>0$, the trajectories of $U(t)$ are
non-decreasing and right-continuous. The process $P(U(t))$,
$t\geq0$, is called a doubly stochastic Poisson process (Cox
process) \cite{Grandell1976}.

Within this model, for each $t$ the distribution of the random
variable $P(U(t))$ is mixed Poisson. For vividness, consider the
case where in the model under consideration the parameter $t$ is
discrete: $U(t)=U(n)=U_n$, $n\in\mathbb{N}$, where $\{U_n\}_{n\ge1}$
is an infinitely increasing sequence of nonnegative random variables
such that $U_{n+1}(\omega)\ge U_{n}(\omega)$ for any
$\omega\in\Omega$, $n\ge1$. Here the asymptotics $n\to\infty$ may be
interpreted as that the intensity of the flow of informative events
is (infinitely) large.

From the assumptions formulated above it follows that the random
variable $U_n$ is independent of the standard Poisson process
$P(t)$, $t\ge0$. For each natural $n$ let $N_n=P(U_n)$, $n\ge1$. It
is obvious that the random variable $N_n$ so defined has the mixed
Poisson distribution
$$
{\sf P}(N_n=k)={\sf
P}\big(P(U_n)=k\big)=\int_{0}^{\infty}e^{-nz}\frac{(nz)^k}{k!}d{\sf
P}(U_n<z)\ \ \ k=0,1,\ldots
$$

Let $X_1,X_2,\ldots $ be independent identically distributed random
variables with the common distribution function $F(x)={\sf
P}(X_i<x)$, $x\in\mathbb{R}$, $i\ge1$. Denote $\mbox{\rm
lext}(F)=\inf\{x:\,F(x)>0\}$. Assume that for each $k\in\mathbb{N}$
the random variable $N_k$ is independent of the sequence
$X_1,X_2,\ldots $ In the book \cite{KorolevSokolov2008} the
following statement was proved.

\smallskip

{\sc Lemma 13} \cite{KorolevSokolov2008}. {\it Assume that there
exist an infinitely increasing sequence of positive numbers
$\{d_k\}_{k\ge1}$ and a nonnegative random variable $U$ such that
$$
\frac{U_k}{d_k}\Longrightarrow U\ \ \ (k\to\infty).
$$
Also assume that $\mbox{\rm lext}(F)>-\infty$ and the distribution
function $A_F(x)=F\bigl( \mbox{\rm lext}(F)-x^{-1}\bigr)$ satisfies
the condition$:$ there exists a positive number $\delta'$ such that
for any $x>0$
$$
\lim_{y\to\-\infty}\frac{A_F(yx)}{A_F(y)}=x^{-\delta'}.\eqno(37)
$$
Then there exist sequences of numbers $a_k$ and $b_k$ such that}
$$
{\sf P}\big(\min_{1\le j\le N_k}X_j-a_k<b_kx\big) \Longrightarrow
\bigg[1-\int_{0}^{\infty}e^{-u x^{\delta'}}d{\sf P}
(U<u)\bigg]\mathbf{1}(x\ge 0)\ \ \ (k\to\infty).
$$
{\it Moreover, the numbers $a_k$ and $b_k$ can be defined as}
$$
a_k={\mbox{\rm lext}}(F),\ \ \ b_k=\sup\big\{x:\ F(x)\le
d_k^{-1}\big\}-{\mbox{\rm lext}}(F),\ \ \ k\ge1.\eqno(38)
$$

\smallskip

From representation (36) and lemma 13 we obtain the following
result.

\smallskip

{\sc Theorem 6.} {\it Let $\delta\in(0,1)$. For the existence of
numbers $a_k\in\r$ and $b_k>0$ such that
$$
\frac{1}{b_k}\Big(\min_{1\le j\le N_k}X_j-a_k\Big)\Longrightarrow
M_{\delta}\ \ \ \ (k\to\infty),
$$
it is sufficient that$:$

\smallskip

\noindent $(i)$ there exists a $\delta'\in(\delta,1]$ such that the
distribution function $F$ belongs to the domain of $\min$-attraction
of the Weibull distribution with some shape parameter
$\delta'\in(0,1]$, that is, $\mbox{\rm lext}(F)>-\infty$ and
condition $(37)$ holds$;$

\smallskip

\noindent $(ii)$ there exists an infinitely increasing sequence
$\{d_k\}_{k\ge1}$ such that}
$$
\frac{U_k}{d_k}\Longrightarrow
S_{\delta',1}^{-\delta'}K_{\delta/\delta'}^{\delta'/\delta}\ \ \ \
(k\to\infty).\eqno(39)
$$

\smallskip

\noindent{\it Moreover, the numbers $a_k$ and $b_k$ can be defined
in accordance with} (38).

\smallskip

{\sc Proof.} The desired result is a direct consequence of lemma 13
and representation (11) with the account of the relation
$K_{\delta/\delta'}^{-1}\eqd K_{\delta/\delta'}$ implied by $(27)$.

\smallskip

So, the randomness of the sample size can make the tails of the
limit distribution considerably more heavy than this is so in the
``classical'' case. For example, let the distribution of the sample
element $X_1$ belong to the domain of min-attraction of the
exponential law, that is, condition (37) holds with $\delta'=1$, but
the sample size is random and has the form $N_k=P(U_k)$ and for some
$d_k$ condition (39) holds with some $\delta\in(0,\delta')$. Then
the actual Mittag-Leffler limit distribution of the minimum order
statistic has power-type decreasing tails unlike the ``originally
assumed'' exponentially decreasing tails.

\subsection{Convergence of the distributions of maximum random sums
to the Mittag-Leffler distribution}

In this section we will demonstrate that the Mittag-Leffler
distribution can be the limit law for maximum or minimum random
sums. The main role here will be played by representations (30) and
(32) of the Mittag-Leffler distribution as a scale mixture of
half-normal distributions. We will show that this distribution can
be limiting for maximum sums of a random number of independent
random variables (maximum random sums), minimum random sums and
absolute values of random sums.

As in Section 6.1, consider independent not necessarily identically
distributed random variables $X_1,X_2,\ldots $ with ${\sf E}X_i=0$
and $0<\sigma^2_i={\sf D}X_i<\infty$, $i\in\N$. In addition to the
notation introduced in Section 6.1, for $n\ge1$ denote $\overline
S^*_n= \max_{1\le i\le n}S_i$, $\underline S^*_n=\min_{1\le i\le
n}S_i$. Assume that the random variables $X_1,X_2,\ldots $ satisfy
the Lindeberg condition (33).

It is well known that under these assumptions ${\sf P}\big(\overline
S^*_n<B_nx\big)\Longrightarrow \Psi(x)$ and ${\sf P}\big(\underline
S^*_n<B_kx\big)\Longrightarrow 1-\Psi(-x)$ as $n\to\infty$ (this is
one of manifestations of the invariance principle).

Let $N_1,N_2,\ldots $ be a sequence of nonnegative random variables
such that for each $n\in\N$ the random variables $N_n,Y_1,Y_2,\ldots
$ are independent. For $n\in\N$ let $S^*_{N_n}=X_1+\ldots +X_{N_n}$,
$\overline S^*_{N_n}= \max_{1\le i\le N_n}S_i$, $\underline
S^*_{N_n}=\min_{1\le i\le N_n}S_i$ (for definiteness assume that
$S_0=\overline S_0=\underline S_0=0$). Let $\{d_n\}_{n\ge1}$ be an
infinitely increasing sequence of positive numbers.

\smallskip

{\sc Lemma 14} \cite{Korolev1994}. {\it Assume that the random
variables $X_1,X_2,\ldots$ and $N_1,N_2,\ldots$ satisfy the
conditions specified above. In particular, let Lindeberg condition
$(33)$ hold. Moreover, let $N_n\pto\infty$ as $n\to\infty$. Then the
distributions of normalized extremal random sums and absolute values
of random sums weakly converge to some distributions, that is, there
exist random variables $Y$, $\overline Y$ and $\underline Y$ such
that
$$
\frac{\overline S^*_{N_n}}{d_n}\Longrightarrow \overline Y,\ \ \ \
\frac{\underline S^*_{N_n}}{d_n}\Longrightarrow \underline Y, \ \ \
\frac{|S^*_{N_n}|}{d_n}\Longrightarrow |Y|
$$
as $n\to\infty$ if and only if there exists a nonnegative random
variable $U$ such that
$$
\frac{B^2_{N_k}}{d_k^{2}}\Longrightarrow U\ \ \ \ (n\to\infty).
$$
Moreover, for $x\in\mathbb{R}$ we have}
$$
{\sf P}\big(\overline Y<x\big)={\sf P}\big(|Y|<x\big)= {\sf
E}\Psi\Big(\frac{x}{\sqrt{U}}\Big),\ \ \ \ {\sf P}\big(\underline
Y<x\big) =1-{\sf E}\Psi\Big(-\frac{x}{\sqrt{U}}\Big).
$$

\smallskip

The {\sc proof} of lemma 14 was given in \cite{Korolev1994}.

\smallskip

Lemma 14 and theorem 4 imply the following statement.

\smallskip

{\sc Theorem 7.} {\it Let $\delta\in(0,1]$. Let $H_{\delta}(x)$ be
the distribution function defined in $(31)$. Assume that the random
variables $X_1,X_2,\ldots$ and $N_1,N_2,\ldots$ satisfy the
conditions specified above. In particular, let Lindeberg condition
$(33)$ hold. Moreover, let $N_n\pto\infty$ as $n\to\infty$. Then, as
$n\to\infty$, the following statements are equivalent}:
$$
\frac{\overline S^*_{N_n}}{d_n}\Longrightarrow M_{\delta}; \ \ \
\frac{\underline S^*_{N_n}}{d_n}\Longrightarrow -M_{\delta};\ \ \
\frac{|S^*_{N_n}|}{d_n}\Longrightarrow M_{\delta};\ \ \ {\sf
P}\big(B^2_{N_n}<d^2_nx\big)\Longrightarrow H_{\delta}(x).
$$

\subsection{Examples of appropriate random indices}

The convergence of the distributions of the normalized indices
$N_n/n$ to a special law is a crucial conditions in all the theorems
presented above concerning the convergence of random sums or
statistics constructed from samples with random sizes to the Linnik
and Mittag-Leffler distributions. For example, the convergence of
the distributions of the normalized indices $N_n/n$ to the
Mittag-Leffler distribution $F_{\delta}^{M}$ is the main condition
in theorem 4. Here we will give two examples of the situation where
this condition can hold. The first example is trivial and is based
on the geometric stability of the Mittag-Leffler distribution. The
second example relies on a useful general construction of
nonnegative integer-valued random variables which, under an
appropriate normalization, converge to a given nonnegative (not
necessarily discrete) random variable, whatever the latter is.
Hence, this construction can be correspondingly modified to give
examples of indices possessing the asymptotic properties required in
other convergence criteria presented in this paper.

\smallskip

{\sc Example 1}. Let $\delta\in(0,1)$ be arbitrary. For every
$n\in\mathbb{N}$ let $V_{1/n}$ be a random variable having the
geometric distribution (4) with $p=\frac1n$ independent of the
sequence $Y_1,Y_2,\ldots$ of independent identically distributed
nonnegative random variables such that
$$
n^{-1/\delta}\sum_{j=1}^{V_{1/n}}Y_j\Longrightarrow
M_{\delta}\eqno(40)
$$
as $n\to\infty$. To provide (40), the distributions of the random
variables $Y_1,Y_2,\ldots$ should belong to the domain of the normal
attraction of the one-sided strictly stable law with characteristic
exponent $\delta$. As $N_n$ for each $n\in\mathbb{N}$ take
$$
N_n=\bigg[n^{1-1/\delta}\sum_{j=1}^{V_{1/n}}Y_j\bigg],
$$
where square brackets denote the integer part. Then
$$
\frac{N_n}{n}=n^{-1/\delta}\sum_{j=1}^{V_{1/n}}Y_j-\frac{1}{n}\bigg\{n^{1-1/\delta}\sum_{j=1}^{V_{1/n}}Y_j\bigg\},\eqno(41)
$$
where curly braces denote the fractional part. Since the second term
on the right-hand side of (41) obviously tends to zero in
probability, from (40) it follows that $N_n/n\Longrightarrow
M_{\delta}$ as $n\to\infty$.

\smallskip

{\sc Example 2.} Consider the construction proposed in Section 6.3.
Assume that the random variable $M_{\delta}$ is independent of the
standard Poisson process $P(t)$, $t\ge0$. For each natural $n$ take
$U_n=nM_{\delta}$. Respectively, let $N_n=P(U_n)=P(nM_{\delta})$,
$n\ge1$. It is obvious that the random variable $N_n$ so defined has
the mixed Poisson distribution
$$
{\sf P}(N_n=k)={\sf
P}\big(P(nM_{\delta})=k\big)=\int_{0}^{\infty}e^{-nz}\frac{(nz)^k}{k!}dF_{\delta}^{M}(z)\
\ \ k=0,1,\ldots
$$
This random variable $N_n$ can be interpreted as the number of
events registered up to time $n$ in the Poisson process with the
stochastic intensity having the Mittag-Leffler density
$f_{\delta}^{M}$.

Denote $A_n(z)={\sf P}(N_n<nz)$, $z\ge0$ ($A_n(z)=0$ for $z<0$). It
is easy to see that $A_n(z)\Longrightarrow F_{\delta}^{M}(z)$ as
$n\to\infty$. Indeed, as is known, if $\Pi(x;\ell)$ is the Poisson
distribution function with the parameter $\ell>0$ and $E(x;c)$ is
the distribution function with a single unit jump at the point
$c\in\r$, then $\Pi(\ell x;\ell)\Longrightarrow E(x;1)$ as
$\ell\to\infty$. Since for $x\in\r$
$$
A_n(x)=\int_{0}^{\infty}\Pi(n x; n z)dF_{\delta}^{M}(z),
$$
then by the Lebesgue dominated convergence theorem, as $n\to\infty$,
we have
$$
A_n(x)\Longrightarrow\int_{0}^{\infty}E(x/z;1)dF_{\delta}^{M}(z)=
\int_{0}^{x}dF_{\delta}^{M}(z)=F_{\delta}^{M}(x),
$$
that is, the random variables $N_n$ defined above satisfy the
condition of theorem 4. Moreover, $N_n\pto\infty$ as $n\to\infty$
since ${\sf P}(M_{\delta}=0)=0$.

As it has already been noted, instead of $M_{\delta}$ any other
positive random variable can be taken to provide the convergence of
the distributions of indices constructed according to the presented
construction to the required distributions.

\section{Two-sided Mittag-Leffler and one-sided Linnik distributions}

\subsection{The symmetric two-sided Mittag-Leffler distribution}

{\sc Definition 1.} Let $\delta\in(0,1)$. By the {\it symmetric
two-sided Mittag-Leffler distribution} with parameter $\delta$ we
will mean the distribution of the random variable $\widetilde
M_{\delta}$ defined by the density
$$
\widetilde
f^{M}_{\delta}(x)=\begin{cases}\frac12f^{M}_{\delta}(-x),&
x<0,\vspace{2mm}\cr \frac12f^{M}_{\delta}(x),&
x\ge0,\end{cases}\eqno(42)
$$
where the Mittag-Leffler density $f^{M}_{\delta}(x)$ was defined in
(2) for $x\in\mathbb{R}^+$.

The distribution function corresponding to the density $\widetilde
f^{M}_{\delta}(x)$ has the form
$$
\widetilde
F^{M}_{\delta}(x)=\begin{cases}\frac12[1-F^{M}_{\delta}(-x)], &
x<0,\vspace{2mm}\cr\frac12[1+F^{M}_{\delta}(x)], & x\ge0.\end{cases}
$$

\smallskip

The two-sided Mittag-Leffler distribution so defined is obviously
symmetric, that is, $-\widetilde M_{\delta}\eqd \widetilde
M_{\delta}$. Furthermore, we can say that $\widetilde M_{\delta}$ is
the {\it randomization symmetrization} of $M_{\delta}$. The
randomization symmetrization of a random variable $Y$ can be
formally defined in the following way.

\smallskip

{\sc Definition 2.} Let $Z$ be a random variable such that ${\sf
P}(Z=-1)={\sf P}(Z=1)=\frac12$. Assume that the random variables $Y$
and $Z$ are independent. The random variable $\widetilde Y=Z\cdot Y$
is called the {\it randomization symmetrization} of $Y$.

\smallskip

It can be easily verified that $\widetilde M_{\delta}\eqd Z\cdot
M_{\delta}$.

\smallskip

It is also easy to verify that the random variable $Y$ is symmetric,
then $\widetilde Y\eqd Z\cdot |Y|$ where the random variables on the
right-hand side are independent and $Z$ is the same as in definition
2.

\smallskip

To describe the symmetric two-sided Mittag-Leffler distribution in
more detail, we will use the following well-known fact (see, e.g.,
\cite{KorolevBeningShorgin2011}). Let $\mathfrak{f}^Y(t)$ be the
characteristic function of the random variable $Y$, and
$\mathfrak{f}^{\widetilde{Y}}(t)$ be the characteristic function of
the randomization symmetrization of the random variable $Y$,
$t\in\r$. Then from definition 2 it follows that
$$
\mathfrak{f}^{\widetilde{Y}}(t)={\textstyle\frac12{\sf
E}e^{-itY}+\frac12{\sf E}e^{itY}}=
{\textstyle\frac12[{\sf E}\cos(tY)-i{\sf E}\sin(tY)]+\frac12[{\sf
E}\cos(tY)+{\sf E}\sin(tY)]}=
$$
$$
={\sf E}\cos(tY)=\Re\mathfrak{f}^Y(t).
$$
That is, the characteristic function of the randomization
symmetrization of any random variable $Y$ coincides with the real
part of the characteristic function of $Y$.

Hence, to find the characteristic function
$\widetilde{\mathfrak{f}}^{M}_{\delta}(t)$ of the symmetric
two-sided Mittag-Leffler distribution, we should calculate the real
part of the characteristic function $\mathfrak{f}^{M}_{\delta}(t)$
of the ordinary (one-sided) Mittag-Leffler distribution.

It is clear that $\mathfrak{f}^{M}_{\delta}(t)=\psi_{\delta}(-it)$,
$t\in\r$ (see (1)). So, we have
$$
\mathfrak{f}^{M}_{\delta}(t)=\psi_{\delta}(-it)=\frac{1}{1+(-it)^{\delta}}=\frac{1}{1+|t|^{\delta}e^{-i\frac{\pi}{2}\delta\,\sign\,t}}=
$$
$$
=\frac{1}{1+|t|^{\delta}[\cos(\frac{\pi}{2}\,\delta\,\sign\,t)-i\sin(\frac{\pi}{2}\delta\,\sign\,t)]}=
\frac{1+|t|^{\delta}\cos(\frac{\pi}{2}\,\delta)+i|t|^{\delta}\sin(\frac{\pi}{2}\,\delta\,\sign\,t)}
{(1+|t|^{\delta}\cos(\frac{\pi}{2}\,\delta\,\sign\,t))^2+
|t|^{2\delta}(\sin(\frac{\pi}{2}\,\delta\,\sign\,t))^2}=
$$
$$
=\frac{1+|t|^{\delta}\cos(\frac{\pi}{2}\,\delta)}{1+|t|^{2\delta}+2|t|^{\delta}\cos(\frac{\pi}{2}\,\delta)}+
i\frac{|t|^{\delta}\sin(\frac{\pi}{2}\,\delta\,\sign\,t)}{1+|t|^{2\delta}+2|t|^{\delta}\cos(\frac{\pi}{2}\,\delta)}.
$$
That is,
$$
\widetilde{\mathfrak{f}}^{M}_{\delta}(t)=\Re\mathfrak{f}^{M}_{\delta}(t)=
\frac{1+|t|^{\delta}\cos(\frac{\pi}{2}\,\delta)}{1+|t|^{2\delta}+2|t|^{\delta}\cos(\frac{\pi}{2}\,\delta)},\
\ \ t\in\r.\eqno(43)
$$
Having compared the right-hand side of (43) with the characteristic
function $\mathfrak{f}^L_{2\delta}(t)$ of the Linnik distribution
(see (5)) we notice that the former differs from the latter by the
presence of the addends of the form
$|t|^{\delta}\cos(\frac{\pi}{2}\,\delta)$ in both the numerator and
denominator. These addends vanish as $\delta\to1$, that is, as both
the symmetric two-sided Mittag-Leffler and Linnik laws turn into the
Laplace distribution.

In general, the above reasoning may be treated as a proof of the
fact that the function on the right-hand side of (43) is a
characteristic function (of a probability distribution). The
properties of this distribution will be considered in more detail in
the next subsection.

\subsection{A normal scale mixture representation for the symmetric two-sided Mittag-Leffler distribution}

{\sc Theorem 8.} {\it The real part of the characteristic function
of any mixed exponential distribution is the characteristic function
of a normal scale mixture.}

\smallskip

{\sc Proof.} Let a random variable $Y$ be represented as $Y\eqd
W_1\cdot U$ where $U$ is a nonnegative random variable independent
of the standard exponential random variable $W_1$. This means that
the distribution of $Y$ is mixed exponential, the mixing
distribution being that of $U^{-1}$. From lemma 10 we obtain the
representation
$$
Y=|X|U\sqrt{2W_1}
$$
with all the random variables on the right-hand side independent.
Now, as it has been demonstrated above,
$$
\Re \mathfrak{f}^Y(t)={\sf E}\exp\{it\widetilde Y\}={\sf
E}\exp\{itZ|X|U\sqrt{2W_1}\}={\sf E}\exp\{itX\sqrt{2W_1U^2}\},\ \ \
t\in\r,\eqno(44)
$$
where all the involved random variables are independent. But
relation (44) means that
$$
{\sf P}(\widetilde
Y<x)=\int_{0}^{\infty}\Phi\Big(\frac{x}{\sqrt{u}}\Big)d{\sf
P}(2W_1U^2<u),\ \ \ x\in\r.
$$
The theorem is proved.

\smallskip

Note that the mixing distributions of the original mixed exponential
distribution and the corresponding normal scale mixture are related
in a simple way: the latter is the mixed exponential distribution
with the mixing law being that of the squared original mixing random
variable.

It is easy to make sure that $X\eqd Z\cdot|X|$, where, as above, $Z$
is a random variable such that ${\sf P}(Z=-1)={\sf P}(Z=1)=\frac12$
and the random variables $X$ and $Z$ are independent. So, from
theorem 8, corollary 2 and (27) we obtain the following result.

\smallskip

{\sc Theorem 9.} {\it The symmetric two-sided Mittag-Leffler
distribution $(42)$ is a normal scale mixture}:
$$
\widetilde M_{\delta}\eqd X\sqrt{2 W_1 K_{\delta}^{2/\delta}}\eqd
X\sqrt{8W_1\Big(\frac{S_{\delta,1}}{S'_{\delta,1}}\Big)^2}.\eqno(45)
$$

\smallskip

We see that in theorem 9 the mixing distribution is the same as in
theorem 4 and was denoted as $H_{\delta}(x)$, see $(31)$, so that
the assertion of theorem 9 can be written as
$$
{\sf P}(\widetilde
M_{\delta}<x)=\int_{0}^{\infty}\Phi\Big(\frac{x}{\sqrt{u}}\Big)dH_{\delta}(u),\
\ \ x\in\r\eqno(46)
$$

To compare the symmetric two-sided Mittag-Leffler distribution with
the Linnik distribution, we will use the representations of these
distributions as scale mixtures of the Laplace distribution. For
$0<\alpha<1$ Denote the distribution function corresponding to
density $p_{\alpha}(x)$ (see (26) by $P_{\alpha}(x)$, $x\ge0$. Let
$\delta\in(0,1)$. Then from theorem 2 it follows that
$$
{\sf
P}(L_{2\delta}<x)=\int_{0}^{\infty}F^{\Lambda}\Big(\frac{x}{\sqrt{u}}\Big)dP_{\delta}(u),\
\ \ x\in\mathbb{R}.\eqno(47)
$$
At the same time, from (22) and theorem 9 we obtain the
representation
$$
{\sf P}(\widetilde
M_{\delta}<x)=\int_{0}^{\infty}F^{\Lambda}\Big(\frac{x}{\sqrt{u}}\Big)dP_{\delta}(\sqrt{2u}),\
\ \ x\in\mathbb{R}.\eqno(48)
$$
Representations (47) and (48) differ only by that the power of the
argument of the mixing distribution in (47) is twice the power of
the argument of the mixing distribution in (48) resulting in that
the exponent of the power-type tail asymptotics of the Linnik
distribution (33) is twice greater than that of the symmetric
two-sided Mittag-Leffler distribution (48).

\subsection{Convergence of the distributions of random sums and
statistics constructed from samples with random sizes to the
symmetric two-sided Mittag-Leffler law}

Following the lines of the reasoning used to prove the results of
Sect. 6.1 and 6.2, here we will present similar results concerning
the convergence of the distributions of random sums and statistics
constructed from samples with random sizes to the symmetric
two-sided Mittag-Leffler law.

Product representations for the random variables with the symmetric
two-sided and Mittag-Leffler distributions obtained in the preceding
sections open the way for the construction in this section of a
random-sum central limit theorem with the symmetric two-sided and
Mittag-Leffler distribution as the limit law. Moreover, as in
Sections 6.1 and 6.2, the results of this section have the ``if and
only if'' character.

We again consider a sequence of independent identically distributed
random variables $X_1,X_2,\ldots$ defined on the probability space
$(\Omega,\, \mathfrak{A},\,{\sf P})$. Assume that ${\sf E}X_1=0$,
$0<\sigma^2={\sf D}X_1<\infty$. We again use the notation
$S^*_n=X_1+\ldots+X_n$, $B_n^2=\sigma_1^2+\ldots+\sigma_n^2$,
$n\in\N$. Let $N_1,N_2,\ldots$ be a sequence of integer-valued
nonnegative random variables defined on the same probability space
so that for each $n\ge1$ the random variable $N_n$ is independent of
the sequence $X_1,X_2,\ldots$ For definiteness, in what follows we
assume that $\sum_{j=1}^0=0$.

Using lemma 11 and theorem 9 we obtain the following statement which
is a criterion (that is, {\it necessary and sufficient} conditions)
of the convergence of the distributions of random sums of
independent identically distributed random variables with {\it
finite} variances to the symmetric two-sided and Mittag-Leffler
distribution.

Let $\{d_n\}_{n\ge1}$ be an infinitely increasing sequence of
positive numbers.

\smallskip

{\sc Theorem 10.} {\it Let $\delta\in(0,1]$. {\it Assume that the
random variables $X_1,X_2,\ldots$ and $N_1,N_2,\ldots$ satisfy the
conditions specified above. In particular, let Lindeberg condition
$(33)$ hold. Moreover, let $N_n\pto\infty$ as $n\to\infty$. Then the
distributions of the normalized random sums $S^*_{N_n}$ converge to
the symmetric two-sided and Mittag-Leffler distribution with
parameter $\delta$, that is,
$$
{\sf P}\Big(\frac{S^*_{N_n}}{d_n}<x\Big) \Longrightarrow \widetilde
F_{\delta}^{M}(x)
$$
as $n\to\infty$, if and only if}
$$
{\sf P}\Big(\frac{B_n^2}{d_n^2}<x\Big)\Longrightarrow H_{\delta}(x)\
\ \ (n\to\infty).
$$
}

\smallskip

Now let random variables $N_1,N_2,\ldots,X_1,X_2,\ldots,$ be defined
on one and the same probability space $(\Omega,{\cal A}, {\sf P})$.
Assume that for each $n\ge 1$ the random variable $N_n$ takes only
natural values and is independent of the sequence $X_1,X_2,\ldots$
Let $T_n=T_n(X_1,\ldots,X_n)$ be a statistic, that is, a measurable
function of $X_1,\ldots,X_n$. Recall that for every $n\ge1$ the
random variable $T_{N_n}$ is defined as
$$
T_{N_n}(\omega)=
T_{N_n(\omega)}\left(X_1(\omega),\ldots,X_{N_n(\omega)}(\omega)\right)
$$
for each $\omega\in\Omega$. We will assume that the statistic
$T_{N_n}$ is asymptotically normal in the sense of (34). Using lemma
12 and theorem 9 we obtain the following criterion (that is, {\it
necessary and sufficient} conditions) for the convergence of the
distributions of statistics, which are suggested to be
asymptotically normal in the traditional sense but are constructed
from samples with random sizes, to the symmetric two-sided
Mittag-Leffler distribution.

\smallskip

{\sc Theorem 11.} {\it Let $\delta\in(0,1]$. Assume that the random
variables $X_1,X_2,\ldots$ and $N_1,N_2,\ldots$ satisfy the
conditions specified above and, moreover, let $N_n\pto\infty$ as
$n\to\infty$. Let the statistic $T_n$ be asymptotically normal in
the sense of $(34)$. Then the distribution of the statistic
$T_{N_n}$ constructed from samples with random sizes $N_n$ converges
to the symmetric two-sided Mittag-Leffler distribution
$F^M_{\delta}(x)$ as $n\to\infty$, that is,
$$
{\sf P}\left(\delta\sqrt{n}\bigl(T_{N_n}-\theta\bigr)<x\right)
\Longrightarrow \widetilde F^M_{\delta}(x),
$$
if and only if
$$
{\sf P}(N_n<nx)\Longrightarrow 1-H_{\delta}\Big(\frac{1}{x}\Big)\ \
\ (n\to\infty).\eqno(49)
$$
}

\smallskip

From the absolute continuity of the distribution function
$H_{\delta}(x)$ it follows that condition (49) can be written as
$$
\sup_{x>0}\Big|{\sf
P}(N_n>nx)-H_{\delta}\Big(\frac{1}{x}\Big)\Big|=0.
$$

\subsection{The one-sided Linnik distribution}

In Section 7.1 we noted that if a random variable $Y$ is symmetric,
that is, $Y\eqd -Y$, and $Z$ is a random variable independent of $Y$
such that ${\sf P}(Z=-1)={\sf P}(Z=1)=\frac12$, then $Y\eqd Z|Y|$.
As this is so, the nonnegative random variable $|Y|$ can be treated
as a ``de-symmetrization'' of $Y$.

Following this logic, the distribution of the random variable
$|L_{\alpha}|$ with $\alpha\in(0,2]$ will be called the one-sided
Linnik distribution.

It is easy to see that
$$
\widehat F^L_{\alpha}(x)\equiv {\sf
P}(|L_{\alpha}|<x)=2F^L_{\alpha}(x)-1, \ \ \ x\ge0.
$$

In \cite{Kozubowski1998} it was shown that for the Linnik
distribution density the following integral representation holds:
$$
f^L_{\alpha}(x)=\frac{\sin(\pi\alpha/2)}{\pi}\int_{0}^{\infty}\frac{y^{\alpha}e^{-y|x|}dy}{1+y^{2\alpha}+2y^{\alpha}\cos(\pi\alpha/2)},
\ \ \ x\in\r.
$$
Hence, the density $\widehat f^L_{\alpha}(x)$ of the one-sided
Linnik law has the form
$$
\widehat
f^L_{\alpha}(x)=\frac{2\sin(\pi\alpha/2)}{\pi}\int_{0}^{\infty}\frac{y^{\alpha}e^{-yx}dy}{1+y^{2\alpha}+2y^{\alpha}\cos(\pi\alpha/2)},
\ \ \ x\ge0.
$$

From corollary 5 we obviously obtain the representation
$$
|L_{\alpha}|=|X|\sqrt{M_{\alpha/2}}.\eqno(50)
$$
Based on this representation, below we will present the conditions
for the one-sided Linnik distribution to be the limit law for
statistics constructed from samples with random sizes, namely, for
maximum random sums and extreme order statistics.

\subsection{Convergence of the distributions of maximum
random sums to the one-sided Linnik law}

In this section we will demonstrate that the one-sided Linnik
distribution can be the limit law for maximum or minimum random
sums. The main role here will be played by representation (50) of
the one-sided Linnik law as a scale mixture of half-normal
distributions. The results of this section are complementary to
those of Section 6.1.

We will use the same notation as in Section 6.4. As in Sections 6.1
and 6.4, assume that the random variables $X_1,X_2,\ldots $ satisfy
the Lindeberg condition (33)

Let $N_1,N_2,\ldots $ be a sequence of nonnegative random variables
such that for each $n\in\N$ the random variables $N_n,X_1,X_2,\ldots
$ are independent. Let $\{d_n\}_{n\ge1}$ be an infinitely increasing
sequence of positive numbers.

Lemmas 11 and 14 together with representation (50) with the account
of the identifiability of scale mixtures of half-normal laws imply
the following statement.

\smallskip

{\sc Theorem 12.} {\it Let $\alpha\in(0,2]$. Assume that the random
variables $X_1,X_2,\ldots$ and $N_1,N_2,\ldots$ satisfy the
conditions specified above. In particular, let Lindeberg condition
$(33)$ hold. Moreover, let $N_n\pto\infty$ as $n\to\infty$. Then, as
$n\to\infty$, the following statements are equivalent}:
$$
\frac{B^2_{N_n}}{d^2_n}\Longrightarrow M_{\alpha/2};\ \ \ \
\frac{S^*_{N_n}}{d_n}\Longrightarrow L_{\alpha}; \ \ \
\frac{\overline S^*_{N_n}}{d_n}\Longrightarrow |L_{\alpha}|; \ \ \
\frac{\underline S^*_{N_n}}{d_n}\Longrightarrow -|L_{\alpha}|;\ \ \
\frac{|S_{N_n}|}{d_n}\Longrightarrow |L_{\alpha}|.
$$

\subsection{Convergence of the distributions of extreme
order statistics to the one-sided Linnik law}

Let $\alpha\in(0,2]$. Since obviously $W_1\eqd|\Lambda|$, from
corollaries 3 and 4 it follows that for any $\gamma\ge1$ we have
$$
|L_{\alpha}|\eqd W_1 Q_{\alpha,2}\eqd W_{\gamma}
T_{1/\gamma}^{1/\gamma} Q_{\alpha,2}\eqd
\frac{2^{1/\gamma}W_{\gamma}}{S_{1/\gamma,1}^{1/\gamma}}\sqrt{\frac{S_{\alpha/2,1}}{S'_{\alpha/2,1}}}.\eqno(51)
$$
We will use the same construction as in Section 6.3. Assume that for
each $n\in\mathbb{N}$ the random variable $U_n$ is independent of
the standard Poisson process $P(t)$, $t\ge0$, and let $N_n=P(U_n)$.
Let $X_1,X_2,\ldots $ be independent identically distributed random
variables with the common distribution function $F(x)={\sf
P}(X_i<x)$, $x\in\mathbb{R}$, $i\ge1$. Denote $\mbox{\rm
lext}(F)=\inf\{x:\,F(x)>0\}$. Assume that for each $n\in\mathbb{N}$
the random variable $N_n$ is independent of the sequence
$X_1,X_2,\ldots $

From representation (51) and lemma 13 we obtain the following
result.

\smallskip

{\sc Theorem 13.} {\it Let $\alpha\in(0,2]$. For the existence of
numbers $a_n\in\r$ and $b_n>0$ such that
$$
\frac{1}{b_n}\Big(\min_{1\le j\le N_n}X_j-a_n\Big)\Longrightarrow
|L_{\alpha}|\ \ \ \ (n\to\infty),
$$
it is sufficient that$:$

\smallskip

\noindent $(i)$ there exists a $\gamma\ge1$ such that the
distribution function $F$ belongs to the domain of $\min$-attraction
of the Weibull distribution with shape parameter $\gamma$, that is,
$\mbox{\rm lext}(F)>-\infty$ and condition $(37)$ holds with
$\delta'=\gamma;$

\smallskip

\noindent $(ii)$ there exists an infinitely increasing sequence
$\{d_n\}_{n\ge1}$ such that}
$$
\frac{U_n}{d_n}\Longrightarrow T_{1/\gamma}^{-1/\gamma}
Q_{\alpha,2}\eqd\frac{S_{1/\gamma,1}^{1/\gamma}}{2^{1/\gamma}}\sqrt{\frac{S_{\alpha/2,1}}{S'_{\alpha/2,1}}}\
\ \ \ (k\to\infty).
$$
{\it where all the random variables on the right-hand side are
independent. Moreover, the numbers $a_n$ and $b_n$ can be defined in
accordance with} (38).

\smallskip

{\sc Proof.} The desired result is a direct consequence of lemma 13
and representation (51) with the account of the relation
$Q_{\alpha,2}^{-1}\eqd Q_{\alpha,2}$ implied by (27).

\renewcommand{\refname}{References}

\end{document}